\documentclass[11pt]{article} 
\usepackage[final]{graphicx}
\usepackage[dvips]{color}
\usepackage{amsmath}
\usepackage{amssymb}
\usepackage{ascmac}
\usepackage{amsthm}
\usepackage{bm}
\usepackage{url}
\usepackage{tabularx}
\usepackage{fancyhdr}
\usepackage[dvipdfmx,svgnames]{xcolor}
\usepackage{tikz}
\usepackage{here}
\usetikzlibrary{calc}
\usepackage{pgfplots}
\usetikzlibrary{positioning}
\def\Bbb{\bf} 

\newcommand\Z{{\Bbb Z}}
\newcommand\Q{{\Bbb Q}}

\newcommand\al{{\alpha}}

\newcommand\g{{\gamma}}  
\newcommand\G{{\Gamma}} 

\newmuskip\pFqskip
\mathchardef\pFcomma=\mathcode`, 

\newcommand*\pFq[5]{%
  \begingroup
  \begingroup\lccode`~=`,
    \lowercase{\endgroup\def~}{\pFcomma\mkern\pFqskip}%
  \mathcode`,=\string"8000
  {}_{#1}F_{#2}\biggl(\genfrac..{0pt}{}{#3}{#4};#5\biggr)%
  \endgroup
}
\newtheorem{lmm}{Lemma}[section]

\newtheorem{rmk}[lmm]{Remark}

\newtheorem{exmp}{Example}
\def\comment#1{ }

\makeatletter
    
    \@addtoreset{equation}{section}
  \makeatother
\allowdisplaybreaks
\begin{document}
\title{Special values of Gauss's hypergeometric series derived from 
Appell's series $F_1$ with closed forms}
\author{Akihito Ebisu
}

\maketitle
\begin{abstract}
In a previous work (\cite{Eb}), 
the author proposed a method employing contiguity relations
to derive hypergeometric series in closed form.
In \cite{Eb}, this method was used to derive 
Gauss's hypergeometric series $_2F_1$ possessing closed forms.
Here, we consider 
the application of this method to Appell's hypergeometetric series $F_1$
and derive several $F_1$ possessing closed forms.
Moreover, analyzing these $F_1$,
we obtain values of $_2F_1$ with no free parameters.
Some of these results provide new examples of algebraic values of $_2F_1$.  

Key Words and Phrases: Gauss's hypergeometric series, algebraic value, 
Appell's hypergeometric series, hypergeometric identity.

2010 Mathematics Subject Classification Numbers: 
Primary 33C05 Secondary 13P15, 33C65, 33F10. 
\end{abstract}
\section{Introduction}
The study of special values of Gauss's hypergeometric series
\begin{gather*}
\pFq{2}{1}{a,b}{c}{x}:=
\sum _{n=0} ^{\infty}\frac{(a,n)(b,n)}{(c,n)(1,n)}x^n,
\label{hgs}
\end{gather*}
where $(\al,n):=\al (\al+1)\cdots (\al +n-1)$,
has a long history.
The oldest and most well-known formula is
\begin{gather}
\pFq{2}{1}{a,b}{c}{1}=
\frac{\G (c)\G(c-a-b)}{\G(c-a)\G(c-b)},
\label{gauss}
\end{gather}
where $\Re (c-a-b)>0$.
The formula (\ref{gauss}), due to Gauss, is called ``Gauss's summation formula''.
Since this formula was derived by Gauss,
many other identities for Gauss's hypergeometric series have been obtained 
by many other people.
For example,  we have
\begin{gather*}
\pFq{2}{1}{2a,2b}{a+b+1/2}{\frac{1}{2}}=
\frac{\sqrt{\pi}\, \G(a+b+1/2)}{\G(a+1/2)\G(b+1/2)}
\label{Kummer}
\end{gather*}
(see formula (50) in Section 2.8 of \cite{Erd1}),
\begin{gather}
\pFq{2}{1}{1/2,-a}{2a+5/2}{\frac{1}{4}}=
\frac{1}{3\cdot2^{2a}}
\frac{\sqrt{\pi}\, \G(2a+5/2)}{\G(a+3/2)^2}
\label{Gosper}
\end{gather}
(see formula (1/4.2) in \cite{Gos} and formula (1,3,3-1)(xvi) in \cite{Eb}).
There are many other known identities for $_2F_1$ similar to the above,
containing one or more free parameters.
These identities can be found by
demonstrating that the corresponding  hypergeometric series possess closed forms.
For instance,
the formula (\ref{Gosper}) is derived by 
showing that $F(a):={}_2F_{1}(1/2, -a; 2a+5/2; 1/4)$ has a closed form, that is, 
that $F(a)$ satisfies the closed-form relation
\begin{gather}
\frac{F(a+1)}{F(a)}=\frac{(2a+5/2,2)}{2^2(a+3/2)^2}.
\label{Gosper2}
\end{gather}
The relation (\ref{Gosper2}) can be obtained by employing
Gosper's algorithm, the W-Z method,
Zeilberger's algorithm (see \cite{Ko} and \cite{PWZ}),
and  the method of contiguity relations, 
which was recently introduced in \cite{Eb}.
Thus, using these methods,
we are able to find numerous identities for $_2F_1$ with one or more free parameters.
Most of the identities that can be derived with these methods
are listed in \cite{Eb}.

There are also many known identities for $_2F_1$ with no free parameters.
The following are some examples:
\begin{gather}
\pFq{2}{1}{1/4,1/2}{3/4}{\frac{80}{81}}=\frac{9}{5}
\label{jz2}
\end{gather}
(see formula (1.6) in \cite{JZ2}),
\begin{gather}
\pFq{2}{1}{1/12,5/12}{1/2}{\frac{1323}{1331}}=
\frac{3\sqrt[4]{11}}{4}
\label{bw}
\end{gather}
(see Theorem 3 in \cite{BW}).
Unfortunately, it seems that
such identities can not be obtained by direct application of the above methods.
Indeed, if we could find (\ref{jz2}) directly with one of the above methods,
then there must exist $p, q, r\in \Q$ satisfying
\begin{gather*}
\frac{F(a+1)}{F(a)}\in \Q(a),
\end{gather*}
where 
\begin{gather*}
F(a):=\pFq{2}{1}{pa+1/4, qa+1/2}{ra+3/4}{\frac{80}{81}}.
\end{gather*}
However, 
no such parameters $p,q,r$ have yet been identified.
The same holds for other identities, including (\ref{bw}).
To obtain these formulae,
other methods have been used, including methods
employing modular forms and elliptic integrals
(see \cite{Ar}, \cite{BG}, \cite{BW}, \cite{JZ1}, \cite{JZ2}, \cite{JZ3}).

As another approach, 
if we could find Appell's hypergeometric series
\begin{gather*}
\pFq{}{1}{\al; \beta _1, \beta_2}{\g}{x,y}
:=\sum _{m=0} ^{\infty}\sum _{n=0} ^{\infty}
\frac{(\al, m+n)(\beta _1,m)(\beta _2,n)}{(\g,m+n)(1,m)(1,n)}x^my^n
\end{gather*}
in closed forms,
then it may be possible to obtain
identities for $_2F_1$ with no free parameters
by analyzing the corresponding closed-form relations for $F_1$.
Below, we consider two examples.

As the First example,
we consider the following
(see Example 1 in Section 3.1 for details).
Using the method of contiguity relations,
which is effective for deriving hypergeometric series in closed form,
we find the closed-form relation
\begin{gather}
\frac{F(a+1)}{F(a)}=\frac{3^{8}}{2^2\, 5^5}\, \frac{(2a+1/2, 2)}{(a+1/2)^2},
\label{ex1closed1}
\end{gather}
where 
\begin{gather*}
F(a):=
\pFq{}{1}
{2a; a+1/2,4a-1}{2a+1/2}{\frac{1}{81},
\frac{1}{6}}.
\end{gather*}
From (\ref{ex1closed1}), 
we obtain 
\begin{gather}
\frac{F(a+n)}{F(a)}=\frac{3^{8n}}{2^{2n}\, 5^{5n}}\, \frac{(2a+1/2, 2n)}{(a+1/2,n)^2},
\label{ex1closed2}
\end{gather}
and
\begin{gather}
\pFq{}{1}
{2a; a+1/2,4a-1}{2a+1/2}{\frac{1}{81},
\frac{1}{6}}
=F(a)=
\frac{3^{8a}}{2^{2a}5^{5a}}\,
\frac{\sqrt{\pi}\,\G(2a+1/2)}{\G(a+1/2)^2}.
\label{ex1gamma}
\end{gather}
Then, from the relation
\begin{gather}
\pFq{}{1}{\al; \beta _1, 0}{\g}{x,y}
=
\pFq{2}{1}{\al; \beta _1}{\g}{x},
\label{reduct}
\end{gather}
which follows from the definition of Appell's hypergeometric series $F_1$,
we have 
\begin{gather}
\pFq{2}{1}{1/2,3/4}{1}{\frac {1}{81}}=
\frac{9}{100}
\frac{\sqrt{2}\, 5^{3/4}\,\G(1/4)^2}{\pi ^{3/2}},
\label{jz3}
\end{gather}
by substituting $a=1/4$ into (\ref{ex1gamma}).

As the second example, we consider the following
(see Example 2 in Section 3.2 for details).
From the method of contiguity relations, 
we are also able to derive the closed-form relation
\begin{gather}
\frac{F(a+1)}{F(a)}=\frac{3^4}{5^4}
\label{ex2closed1}
\end{gather}
and, from this, 
\begin{gather}
{F(a)}=\frac{5^{4n}}{3^{4n}}F(a+n),
\label{ex2closed2}
\end{gather}
where 
\begin{gather}
F(a):=
\pFq{}{1}
{a;2a,1-4a}{a+1/2}{\frac{80}{81},
\frac{16}{15}}.
\label{ex2}
\end{gather}
It is known that, in general, $F_1(\al; \beta _1, \beta _2; \g; x, y)$ 
is absolutely convergent
when $|x|<1$ and $|y|<1$, and divergent when $|x|>1$ or $|y|>1$
(for example, see Section 9.1 in \cite{Ba}).
Therefore, $F(a)$ is not meaningful for a parameter $a$ with unrestricted value.
However, $F(a)$ is meaningful if we restrict $a$ to values satisfying
$a=1/4+n$, where $n$ is any non-negative integer.
Thus, the relation (\ref{ex2closed2}) with $a=1/4$ is meaningful.
Then, investigating the asymptotic behavior of 
the right-hand side of this relation
by taking the limit $n\rightarrow +\infty$,
we are able to deduce (\ref{jz2}).
In this way,
analyzing closed-form relations for $F_1$, we can derive values of $_2F_1$
with no free parameters.

In this article, 
making use of the method of contiguity relations,
we derive several $F_1$ in closed forms.
These are listed in Table \ref{closed_form}.
Moreover,  
in the cases that these $F_1$ are convergent,
we evaluate them with (\ref{ex1gamma}).
These hypergeometric identities for $F_1$
are tabulated in 
Table \ref{identities}.
In addition,
analyzing the closed-form relations for $F_1$ listed in Table \ref{closed_form},
we derive values of $_2F_1$ with no free parameters.
These values are listed in Table \ref{value1}.
In Tables \ref{value2} and \ref{value3},
we present several complicated identities obtained by
applying algebraic transformations to the identities listed in Table \ref{value1}.
As seen from these,
we are able to obtain several new identities using our approach.
In particular, (B$''$.3), (C$''$.1), (C$''$.4) and (C$'''$.1)
provide new examples of algebraic values of Gauss's hypergeometric series. 
However, it seems that we are not able to derive 
the beautiful formula (\ref{bw}) obtained by Beukers and Wolfart
with our approach.

\begin{rmk}
In \cite{Si}, Siegel posed the problem of determining the nature of the following set:
\begin{gather*}
E(a,b,c):=\{x\in \bar{\Q};\ {}_2F_1(a,b;c;x)\in \bar{\Q} \}.
\end{gather*} 
This set is a so-called ``exceptional set''.
Since Wolfart's celebrated work \cite{Wo} investigating $E(a,b,c)$,
much progress has been made in the study of exceptional sets.
In those studies, 
Gauss's hypergeometric series $_2F_1(a,b;c;x)$ corresponding to 
the arithmetic triangular group
and satisfying the relations
\begin{gather*}
c<1,\ 0<a<c,\ 0<b<c,\ |1-c|+|a-b|+|c-a-b|<1,\\
\frac{1}{|1-c|}, \frac{1}{|a-b|}, \frac{1}{|c-a-b|}\in \Z _{>0}
\end{gather*}
is a focus of investigation.
The formulae 
{\rm (B$''$.3), (C$''$.1), (C$''$.4)} and {\rm (C$'''$.1)} provide
new examples of elements of these exceptional sets.
The formulae {\rm (B$''$.3), (C$''$.1), (C$''$.4)} and {\rm (C$'''$.1)}  
arise from Gauss's hypergeometric equations corresponding
to the $(3,6,6)$, $(2,5,10)$, $(2,5,10)$ and $(2,3,10)$ triangular groups, respectively.
\end{rmk}

\section{The method of contiguity relations}
In this section,
we review the method of contiguity relations, 
which was introduced in \cite{Eb}.
Using this method,
we obtain the examples of Appell's hypergeometric series $F_1$ 
possessing closed forms listed in Table \ref{closed_form}.

First, for simplicity, we write the parameters $(\al; \beta _1, \beta _2; \g)$
as ${\bm \al}$.
Then, Appell's hypergeometric series 
\begin{gather*}
\pFq{}{1}{\al; \beta _1, \beta _2}{\g}{x,y}
\end{gather*}
is expressed as $F_1({\bm \al}; x,y)$.
We also define 
\begin{gather*}
{\bm e_{10}}:=(1;1,0;1),\
{\bm e_{01}}:=(1;0,1;1),\
{\bm k}:=(k; l _1, l _2; m)\in \Z^4.
\end{gather*}

Now, we review the method of contiguity relations.
It is known that for a given quadruple of integers ${\bm k} \in\Z^4$,
there exists a unique triple of rational functions 
$(Q_{10}, Q_{01}, Q_{00}) \in (\Q (\al; \beta _1, \beta _2; \g, x, y))^3$ 
satisfying 
\begin{align}
\begin{split}
F_1({\bm \al}+{\bm k};x,y)
=Q_{10}\,F_1({\bm \al}+{\bm e_{10}};x,y)
+Q_{01}\,F_1({\bm \al}+{\bm e_{01}};x,y)
+Q_{00}\,F_1({\bm \al};x,y).
\end{split}
\label{four_term}
\end{align}
The relation (\ref{four_term}) is called the ``contiguity relation'' for $F_1$
(or the ``four-term relation'' for $F_1$). 
Note that 
it is possible to compute 
$Q_{10}$, $Q_{01}$ and $Q_{00}$ 
exactly by using the method introduced in \cite{Ta1}
(see also \cite{Ta2}).
Next, we define 
\begin{gather*}
Q_{ij} ^{(n)}:=Q_{ij} |_{{\bm \al} \rightarrow
{\bm \al}+n{\bm k}},
\end{gather*}
where each $Q_{ij} ^{(n)}$ is a rational expression in $n$ 
whose coefficients belong to $\Z[\al,$ $\beta _1$,$ \beta _2$, $\g$, $x$, $y]$.
Now, let the sextuple $(\al; \beta _1, \beta _2; \g, x, y)$ satisfy 
\begin{gather}
\begin{cases}
Q_{10}^{(n)}\equiv 0,\\
Q_{01}^{(n)}\equiv 0\quad
\label{case}
\end{cases}
 {\text {for every integer $n$}}.
\end{gather}
Such a sextuple can be found 
by solving the polynomial system in $(\al, \beta _1, \beta _2, \g, x, y)$ 
arising from (\ref{case}).
Then,
from (\ref{four_term}),
we find that the relation
\begin{align}
F_1({\bm \al}+(n+1){\bm k};x,y)
=Q_{00} ^{(n)}\,F_1({\bm \al}+n{\bm k};x,y)
\label{two_term_f1}
\end{align}
holds for such a sextuple.
The relation (\ref{two_term_f1}) implies that
$F(n):=F_1({\bm \al}+n{\bm k};x,y)$ has a closed form.
Thus, we are able to find $F_1$ in closed form.
The above method is called the ``method of contiguity relations''.

As an example,
we consider the case ${\bm k}=(2,1,4,2)$.
Applying the method of contiguity relations to this case,
we find that
\begin{gather*}
(\al, \beta _1, \beta _2, \g, x, y)
=\left(2a, a+\frac{1}{2}, 4a-1, 2a+\frac{1}{2}, \frac{1}{81}, \frac{1}{6}\right)
\end{gather*}
satisfies (\ref{case}).
For this sextuple, $Q _{00}$ becomes
\begin{gather*}
\frac{3^{8}}{2^2\, 5^5}\, \frac{(2a+1/2, 2)}{(a+1/2)^2}.
\end{gather*}
Hence, we have the following relation:
\begin{align*}
\dfrac{F_1({\bm \al}+(n+1){\bm k};x,y)}{F_1({\bm \al}+n{\bm k};x,y)}
&=
\dfrac{
\pFq{}{1}{2a+2n+2; a+n+3/2,4a+4n+3}{2a+2n+5/2}{\frac{1}{81},\frac{1}{6}}
}
{
\pFq{}{1}{2a+2n; a+n+1/2,4a+4n-1}{2a+2n+1/2}{\frac{1}{81},\frac{1}{6}}
}
\\
&=
\frac{3^{8}}{2^2\, 5^5}\, \frac{(2a+2n+1/2, 2)}{(a+n+1/2)^2}=Q _{00} ^{(n)}.
\end{align*}
Thus, with
\begin{gather*}
F(a):=
\pFq{}{1}
{2a; a+1/2,4a-1}{2a+1/2}{\frac{1}{81},
\frac{1}{6}},
\end{gather*}
we see that $F(a)$ has a closed form, and $F(a)$ satisfies the closed-form relation
\begin{gather*}
\frac{F(a+1)}{F(a)}=
\frac{3^{8}}{2^2\, 5^5}\, \frac{(2a+1/2, 2)}{(a+1/2)^2}.
\end{gather*}

As seen in the above example,
if $\bm k \in \Z^4$ is given,
then we can obtain $F_1$ in closed form by using the method of contiguity relations.
Such examples are tabulated in  Table \ref{closed_form}.

\begin{table}[htb]
\caption{
Some examples of  $F_1$ possessing closed forms. 
} 
\label{closed_form}
\scalebox{1.00}{
\begin{tabular}{llll}
\hline
No. & ${\bm k}$ & $F(a)$ & $F(a+1)/F(a)$ \\ 
\hline
(A.1)
&
$(2,1,4,2)$
&
\shortstack{{} \\
${\displaystyle
\pFq{}{1}
{2a; a+1/2,4a-1}{2a+1/2}{\frac{1}{81},
\frac{1}{6}}
}$}
&
${\displaystyle
\frac{3^{8}}{2^2\, 5^5}\, \frac{(2a+1/2, 2)}{(a+1/2)^2}
}$
\\
(A.2)
&
$(2,1,4,5)$
&
\shortstack{{} \\
${\displaystyle
\pFq{}{1}
{2a;a+1/2,4a-1}{5a}{\frac{80}{81},
\frac{5}{6}}
}$}
&
${\displaystyle
\frac{3^8}{2^2\, 5^5}\,
\frac{(5a,5)}{(a+1/2)^2(3a,3)}
}$
\\
(A.3)
&
$(1,2,-4,1)$
&
\shortstack{{} \\
${\displaystyle
\pFq{}{1}
{a;2a,1-4a}{a+1/2}{\frac{80}{81},
\frac{16}{15}}
}$}
&
${\displaystyle
\frac{3^4}{5^4}
}$
\\
(B.1)
&
$(2,-3,4,2)$
&
\shortstack{{} \\
${\displaystyle
\pFq{}{1}
{2a;1-3a,4a-1}{2a+1/2}{\frac{-1}{80},\frac{5}{32}}
}$}
&
${\displaystyle
\frac{2^6}{5^3}\, \frac{(2a+1/2,2)}{(a+1/2)^2}
}$
\\
(B.2)
&
$(1,3,0,5)$
&
\shortstack{{} \\
${\displaystyle
\pFq{}{1}{a;3a-1/2,1/2}{5a}{\frac{-25}{2}, \frac{5}{32}}
}$}
&
${\displaystyle
\frac{2^3}{5^5}\,
\frac{a(5a,5)}{(2a,2)^2(2a+1/2,2)}
}$
\\
(B.3)
&
$(2,-3,4,1)$
&
\shortstack{{} \\
${\displaystyle
\pFq{}{1}{2a;1-3a,4a-1}{a+1/2}
{\frac{81}{80},\frac{27}{32}}
}$}
&
${\displaystyle
-\frac{2^6}{5^3}
}$
\\
(B.4)
&
$(0,-1,3,4)$
&
\shortstack{{} \\
${\displaystyle
\pFq{}{1}{1/2;-a,3a+5/2}{4a+9/2}{\frac{27}{32},\frac{5}{6}}
}$}
&
${\displaystyle
\frac{(4a+9/2,4)}{2^6\,(a+3/2)^2(2a+5/2,2)}
}$
\\
(C.1)
&
$(2,1,5,3)$
&
\shortstack{{} \\
${\displaystyle
\pFq{}{1}{2a;a+1/2,5a-3/2}{3a}{\frac{3}{128},\frac{3}{8}}
}$}
&
${\displaystyle
\frac{2^{16}}{3^35^5}\,\frac{(3a,3)}{(a+1/2)(2a,2)}
}$
\\
(C.2)
&
$(3,-1,5,2)$
&
\shortstack{{} \\
${\displaystyle
\pFq{}{1}{3a;1-a,5a-3/2}{2a+1/2}{\frac{3}{128},\frac{1}{16}}
}$}
&
${\displaystyle
\frac{2^{15}}{5^5}\,\frac{a(2a+1/2,2)}{(3a,3)}
}$
\\
(C.3)
&
$(2,1,5,5)$
&
\shortstack{{} \\
${\displaystyle
\pFq{}{1}{2a;a+1/2,5a-3/2}{5a}{\frac{125}{128},\frac{5}{8}}
}$}
&
${\displaystyle
\frac{2^{18}}{3^35^5}\,\frac{a^2(5a,5)}{(2a,2)^2(3a,3)}
}$
\\
(C.4)
&
$(1,-3,5,0)$
&
\shortstack{{} \\
${\displaystyle
\pFq{}{1}{a;1-3a,5a-3/2}{1/2}{\frac{125}{128},\frac{25}{16}}
}$}
&
${\displaystyle
\frac{2}{3^3}
}$
\\
(D.1)
&
$(2,-3,5,3)$
&
\shortstack{{} \\
${\displaystyle
\pFq{}{1}{2a;1-3a,5a-3/2}{3a}{\frac{-3}{125},\frac{9}{25}}
}$}
&
${\displaystyle
\frac{2^25}{3^3}\,
\frac{(3a,3)}{(a+1/2)(2a,2)}
}$
\\
(D.2)
&
$(1,1,3,-1)$
&
\shortstack{{} \\
${\displaystyle
\pFq{}{1}{a;a+1/2,3a-1/2}{1-a}{\frac{16}{25},16}
}$}
&
${\displaystyle
\frac{5^2}{3^6}
}$
\\
(D.3)
&
$(5,3,2,4)$
&
\shortstack{{} \\
${\displaystyle
\pFq{}{1}{5a;3a-1/2,2a}{4a+1/2}{\frac{1}{25},\frac{16}{25}}
}$}
&
${\displaystyle
\frac{5^{10}}{2^63^6}\,\frac{a(4a+1/2,4)}{(5a,5)}
}$
\\
(E.1)
&
$(0,2,3,3)$
&
\shortstack{{} \\
${\displaystyle
\pFq{}{1}{1/2;2a,3a-1/2}{3a+1/2}{\frac{1}{5},-\frac{4}{5}}
}$}
&
${\displaystyle
\frac{(a+1/6)(a+5/6)}{(a+1/3)(a+2/3)}
}$
\\
\hline
\end{tabular}
}
\end{table}

\section{Special values of $_2F_1$ derived from $F_1$ with closed forms}
In the previous section,
we presented examples of $F_1$ possessing closed forms.
These are listed in Table \ref{closed_form}.
In this section, from among these $F_1$,
we evaluate those that are convergent.
The derived hypergeometric identities for $F_1$
are listed in Table \ref{identities}. 
In addition, analyzing the closed-form relations for $F_1$
listed in Table \ref{closed_form},
we obtain values of $_2F_1$ with no free parameters.
These values are presented in Table \ref{value1}.

\subsection
{Hypergeometric identities for $F_1$ obtained from Table \ref{closed_form}}
In this subsection,
we derive hypergeometric identities 
from the closed-form relations given in Table \ref{closed_form}.

\begin{exmp}{\rm
As seen in Table \ref{closed_form},
setting ${\bm k}=(2,1,4,2)$,
we find that
\begin{gather*}
F(a):=
\pFq{}{1}
{2a;a+1/2,4a-1}{2a+1/2}{\frac{1}{81},
\frac{1}{6}}
\end{gather*}
possesses a closed form, and it satisfies the closed-form relation
\begin{gather}
\frac{F(a+1)}{F(a)}=
\frac{3^{8}}{2^2\, 5^5}\, \frac{(2a+1/2, 2)}{(a+1/2)^2}.
\label{gamma_ex}
\end{gather}
This implies 
\begin{gather}
\frac{F(a+n)}{F(a)}=
\frac{3^{8n}}{2^{2n}\, 5^{5n}}\, \frac{(2a+1/2, 2n)}{(a+1/2,n)^2}.
\label{gamma_ex1}
\end{gather}
Then, noting that
\begin{gather*}
\pFq{}{1}{0; \beta _1, \beta _2}{\g}{x,y}=1,
\end{gather*}
and, substituting $a=0$ into (\ref{gamma_ex1}),
we obtain
\begin{gather}
F(n)=G(n)
\label{gamma_ex2}
\end{gather}
for any integer $n$, where
\begin{gather*}
G(n)
:=
\frac{3^{8n}}{2^{2n}\, 5^{5n}}\, \frac{(1/2, 2n)}{(1/2,n)^2}
=
\frac{3^{8n}}{2^{2n}\, 5^{5n}}\, \frac{\G(1/2)\G(2n+1/2)}{\G(n+1/2)^2}.
\end{gather*}
Now, we show that the identity (\ref{gamma_ex2}) holds 
for any complex number $a$, except for $a=-1/4, -3/4, -5/4, \ldots$.
For this purpose, we use the following lemma proved by Carlson
(see Section 5.3 in \cite{Ba}).
\begin{lmm}
(Carlson's theorem)
We assume that $f(a)$ and $g(a)$ are regular and of the form $O(e^{k|a|})$, where $k<\pi$, 
for $\Re (a) \geq 0$,
and that $f(a)=g(a)$ for $a=0,1,2,\ldots$. 
Then, we have $f(a)=g(a)$ on $\{a\mid \,\Re (a) \geq 0\}$.
\end{lmm}
\noindent
Obviously, both sides of (\ref{gamma_ex2}) are regular for $\Re (a) \geq 0$.
Also, we have the following: 
\begin{align*}
&|F(a)|\\
&\leq 
\sum _{m,n=0} ^{\infty}
\frac{|2a||2a+1|\cdots |2a+m+n-1|(|a|+1/2,m)(4|a|+1,n)}
{|2a+1/2||2a+3/2|\cdots|2a-1/2+m+n|(1,m)(1,n)}
\left(\frac{1}{81}\right)^m
\left(\frac{1}{6}\right)^n\\
&\leq
\sum _{m,n=0} ^{\infty}
\frac{(|a|+1/2,m)(4|a|+1,n)}{(1,m)(1,n)}
\left(\frac{1}{81}\right)^m
\left(\frac{1}{6}\right)^n\\
&=
\sum _{m=0} ^{\infty}
\frac{(|a|+1/2,m)}{(1,m)}
\left(\frac{1}{81}\right)^m
\sum _{n=0} ^{\infty}
\frac{(4|a|+1,n)}{(1,n)}
\left(\frac{1}{6}\right)^n\\
&=
\left(1-\frac{1}{81}\right)^{-|a|-1/2}
\left(1-\frac{1}{6}\right)^{-4|a|-1}\\
&=
O\left(\exp\left(|a|\log\left(\frac{3^8}{5^5}\right)\right)\right).
\end{align*}
Note that we can also compute the asymptotic behavior of $F(a)$ exactly 
by making use of Laplace's method
(also, see Section 3 in \cite{Iw}).
In any case, we see that $|F(a)|$ is of the form $O(e^{|a|})$,
whereas it is easily demonstrated using Stirling's formula
that $|G(a)|$ is of the form $O(e^{|a|})$.
Of course, we know that from (\ref{gamma_ex2}),
\begin{gather}
F(a)=G(a)
\label{gamma_ex3}
\end{gather}
holds for any non-negative integer $a$.
Therefore, it follows from Carlson's theorem 
that the identity (\ref{gamma_ex3}) 
is valid for $\Re(a) \geq 0$.
Also, by analytic continuation,
we find that (\ref{gamma_ex3}) holds for  
any complex number $a$, except for $a=-1/4, -3/4, -5/4, \ldots$.
Thus, we derive (\ref{ex1gamma}), which appears as (A$'$.1) in Table \ref{identities}.
}\end{exmp}
In this way,
it is possible to obtain hypergeometric identities for $F_1$
from the closed-form relations in Table \ref{closed_form}.
These are listed in Table \ref{identities}. 

\begin{table}[htb]
\caption{
Hypergeometric identities obtained from Table \ref{closed_form}.
} 
\label{identities}
\scalebox{1.00}{
\begin{tabular}{llll}
\hline
No. &  Hypergeometric identity  \\ 
\hline
(A$'$.1)
&
\shortstack{{} \\
${\displaystyle
\pFq{}{1}
{2a;a+1/2,4a-1}{2a+1/2}{\frac{1}{81},
\frac{1}{6}}
=
\frac{3^{8a}}{2^{2a}5^{5a}}\,
\frac{\sqrt{\pi}\,\G(2a+1/2)}{\G(a+1/2)^2}.
}$}
\\
(A$'$.2)
&
\shortstack{{} \\
${\displaystyle
\pFq{}{1}
{2a;a+1/2,4a-1}{5a}{\frac{80}{81},
\frac{5}{6}}
=
\frac{3^{8a}}{2^{2a-1}5^{5a}}\,
\frac{\pi\,\G(5a)}{\G(a+1/2)^2\G(3a)}.
}$}
\\
(B$'$.1)
&
\shortstack{{} \\
${\displaystyle
\pFq{}{1}
{2a;1-3a,4a-1}{2a+1/2}{\frac{-1}{80},\frac{5}{32}}
=\frac{2^{6a}}{5^{3a}}\,
\frac{\sqrt{\pi}\,\G(2a+1/2)}{\G(a+1/2)^2}.
}$}
\\
(B$'$.4)
&
\shortstack{{} \\
${\displaystyle
\pFq{}{1}{1/2;-a,3a+5/2}{4a+9/2}{\frac{27}{32},\frac{5}{6}}
=
\frac{1}{5\sqrt{3}\, 2^{6a+3/2}}\,\frac{\pi\,\G(4a+9/2)}{\G(a+3/2)^2\G(2a+5/2)}.
}$}
\\
(C$'$.1)
&
\shortstack{{} \\
${\displaystyle
\pFq{}{1}{2a;a+1/2,5a-3/2}{3a}{\frac{3}{128},\frac{3}{8}}
=\frac{2^{16a}}{3^{3a}5^{5a}}\,
\frac{\sqrt{\pi}\, \G(3a)}{\G(a+1/2)\G(2a)}.
}$}
\\
(C$'$.2)
&
\shortstack{{} \\
${\displaystyle
\pFq{}{1}{3a;1-a,5a-3/2}{2a+1/2}{\frac{3}{128},\frac{1}{16}}
=\frac{2^{15a}}{3\cdot5^{5a}}\,\frac{\G(a)\G(2a+1/2)}{\sqrt{\pi}\, \G(3a)}.
}$}
\\
(C$'$.3)
&
\shortstack{{} \\
${\displaystyle
\pFq{}{1}{2a;a+1/2,5a-3/2}{5a}{\frac{125}{128},\frac{5}{8}}
=\frac{2^{18a-1}}{3^{3a}5^{5a}}\,\frac{\G(a)^2\G(5a)}{\G(2a)^2\G(3a)}.
}$}
\\
(D$'$.1)
&
\shortstack{{} \\
${\displaystyle
\pFq{}{1}{2a;1-3a,5a-3/2}{3a}{\frac{-3}{125},\frac{9}{25}}
=\frac{2^{2a}5^a}{3^{3a}}\,
\frac{\sqrt{\pi}\,\G(3a)}{\G(a+1/2)\G(2a)}.
}$}
\\
(D$'$.3)
&
\shortstack{{} \\
${\displaystyle
\pFq{}{1}{5a;3a-1/2,2a}{4a+1/2}{\frac{1}{25},\frac{16}{25}}
=\frac{5^{10a-1}}{2^{6a}3^{6a}}\,\frac{\G(a)\G(4a+1/2)}{\sqrt{\pi}\, \G(5a)}.
}$}
&
\\
(E$'$.1)
&
\shortstack{{} \\
${\displaystyle
\pFq{}{1}{1/2;2a,3a-1/2}{3a+1/2}{\frac{1}{5},-\frac{4}{5}}
=
\sqrt{\frac{3}{5}}\,
\frac{\G(a+1/6)\G(a+5/6)}{\G(a+1/3)\G(a+2/3)}.
}$}
\\
\hline
\end{tabular}
}
\end{table}

\subsection{Special values of $_2F_1$ derived from closed-form relations for $F_1$}
In this subsection,
we derive identities for $_2F_1$  with no free parameters,
using the closed-form relations listed in Table \ref{closed_form}
and the hypergeometric identities for $F_1$ listed in Table \ref{identities}. 

For example,
let us substitute $a=1/4$ into (A$'$.1).
Then, recalling the reduction formula (\ref{reduct}), 
we obtain (\ref{jz3}).
This appears as formula (A$''$.1) in Table \ref{value1}. 
In a similar way,
reducing the hypergeometric identities for $F_1$ 
appearing in Table \ref{identities}
to identities for $_2F_1$,
we derive 
(A$''$.2), (B$''$.1), (B$''$.4), (C$''$.1), (C$''$.2), (C$''$.3),
(D$''$.1), (D$''$.3) and (E$''$.1).

We now derive the remaining identities 
that can be obtained from the relations and identities given 
in Tables \ref{closed_form} and \ref{identities}.
\begin{exmp}{\rm
From Table \ref{closed_form}, choosing ${\bm k}=(1,2,-4,1)$, we find that
\begin{gather*}
F(a):=
\pFq{}{1}
{a;2a,1-4a}{a+1/2}{\frac{80}{81},
\frac{16}{15}}
\end{gather*}
has a closed form, and it satisfies the following closed-form relation:
\begin{gather}
\frac{F(a+1)}{F(a)}=\frac{3^4}{5^4}.
\label{ex2_closed1}
\end{gather}
This implies that
\begin{gather}
{F(a)}=\frac{5^{4n}}{3^{4n}}{F(a+n)}
\label{ex2_closed2}
\end{gather}
holds for any non-negative integer $n$.
Although (\ref{ex2_closed2}) (and (\ref{ex2_closed1})) is valid
by virtue of analytic continuation,
$F(a)$ and $F(a+n)$,
which are regarded as infinite double series expressions, are meaningless. 
For this reason, we carry out a reduction of each of these
to a finite sum of a single series expression that is meaningful.
This is done, for example, by substituting $a=1/4$ into (\ref{ex2_closed2}).
Then, it becomes 
\begin{gather}
\pFq{2}{1}{1/4, 1/2}{3/4}{\frac{80}{81}}
=
\frac{5^{4n}}{3^{4n}}
\pFq{}{1}{1/4+n; 1/2+2n, -4n}{3/4+n}{\frac{80}{81}, \frac{16}{15}}.
\label{ex2_eq}
\end{gather}
Next, we evaluate the left-hand side of (\ref{ex2_eq})
by determining the asymptotic behavior of the right-hand side of (\ref{ex2_eq})
in the limit $n\rightarrow +\infty$. 
Because $F_1(\al; \beta _1, \beta _2; \g; x,y)$ has the integral representation
\begin{gather*}
\pFq{}{1}{\al; \beta _1, \beta _2}{\g}{x,y}
=\frac{\G(\g)}{\G(\al)\G(\g-\al)}
\int _{0} ^{1} t^{\al-1}(1-t)^{\g-\al-1}(1-xt)^{-\beta _1}(1-yt)^{-\beta _2}dt,
\end{gather*}
the right-hand side of (\ref{ex2_eq}) can be expressed as $A\cdot B$, where
\begin{align*}
A:=\frac{5^{4n}}{3^{4n}}
\frac{\G(3/4+n)}{\G(1/4+n)\G(1/2)},\ 
B:=\int _{0} ^{1}
g(t)e ^{nh(t)}dt,
\end{align*}
and, here,
\begin{gather*}
g(t):=t^{-3/4}(1-t)^{-1/2}\left(1-\frac{80}{81}t\right)^{-1/2},\
h(t):=\log\left(
t\left(1-\frac{80}{81}t\right)^{-2}\left(1-\frac{16}{15}t\right)^4
\right).
\end{gather*}
The function $h(t)$ is plotted in Figure \ref{fig1}.
\begin{figure}[htbp]
 \begin{center}
  \includegraphics[width=80mm]{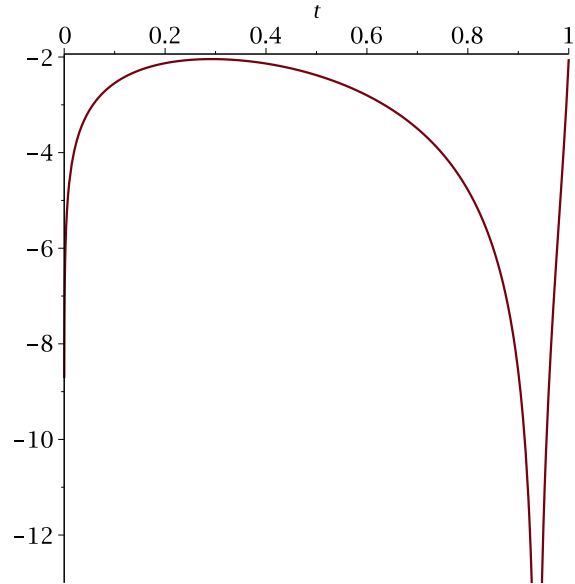}
 \end{center}
 \caption{Graph of $h(t)$.}
 \label{fig1}
\end{figure}
From the relation
\begin{gather*}
h'(t)={\frac {15(256\,{t}^{2}-352\,t+81)}{t \left( 15-16\,t \right) 
 \left( 81-80\,t \right) }},
\end{gather*}
we find that $h'(t)$ becomes zero at 
\begin{gather*}
t_0:=\frac{11}{16}-\frac{\sqrt{10}}{8}.
\end{gather*}
We also obtain the relations
\begin{gather*}
h(t_0)=h(1)=\log\left(\frac{3^4}{5^4}\right).
\end{gather*}
From Figure \ref{fig1}, it can be seen that
the major contribution to the value of $B$ arises 
from the neighborhoods of the points $t=t_0, 1$.
So, using Laplace's method
(see Section 2.4 in \cite{Erd2}),
we find that $B$ takes the form
\begin{gather}
B\sim \frac{9}{5}\left(\frac{3^4}{5^4}\right)^{n}\sqrt{\frac{\pi}{n}}
\label{ex2_asymb}
\end{gather}
in the limit $n \rightarrow +\infty$.
We can also compute the asymptotic behavior of $A$ using Stirling's formula;
we find that
\begin{gather}
A\sim \left(\frac{5^4}{3^4}\right)^{n}\sqrt{\frac{n}{\pi}}
\label{ex2_asyma}
\end{gather}
in the limit $n \rightarrow +\infty$.
The formulae (\ref{ex2_asymb}) and (\ref{ex2_asyma}) yield
\begin{gather*} 
\lim _{n\rightarrow +\infty}A\cdot B=\frac{9}{5}.
\end{gather*}
Thus, we obtain (\ref{jz2}).
This appears as (A$''$.3) in Table \ref{value1}.
}\end{exmp}

Similarly, 
some of the remaining identities can be obtained.
However, among the relations given in Table \ref{closed_form},
there are some cases 
in which it it is difficult to obtain values of $_2F_1$
by direct application of Laplace's method.
For such cases,
we must use connection formulae for $_2F_1$.

For example, 
it is apparently difficult to compute the asymptotic behavior of 
(B.3) with $a=1/3+n$ in the limit $n\rightarrow +\infty$,
by direct use of Laplace's method.
For this reason, we use a connection formula for $_2F_1$
to obtain (B$''$.3) in Table \ref{value1}.
As seen in the beginning of this subsection,
we already know (B$''$.1) and (B$''$.4).
Also, the following is a known connection formula for $_2F_1$:
\begin{gather}
u_1=\frac{\G(c)\G(c-a-b)}{\G(c-a)\G(c-b)}u_2
+\frac{\G(c)\G(a+b-c)}{\G(a)\G(b)}u_6,
\label{connection}
\end{gather}  
where
\begin{align*}
u_1&:=\pFq{2}{1}{a,b}{c}{x},\\
u_2&:=\pFq{2}{1}{a,b}{a+b+1-c}{1-x},\\
u_6&:=(1-x)^{c-a-b}\pFq{2}{1}{c-a,c-b}{c+1-a-b}{1-x}
\end{align*}
(see formula (33) in Section 2.9 in \cite{Erd1}).
Then, substituting $(a,b,c,x)=(1/3$, $2/3,$ 7/6$, 5/32)$ into (\ref{connection}),
we deduce (B$''$.3).
In this way, all the remaining identities are derived.
All of the identities we have been able to derive are presented in Table \ref{value1}.
There, ``type'' refers to the type of Schwarz triangle for 
the Schwarz map of the corresponding hypergeometric equation.
Explicitly, for a given $_2F_1(a,b;c;x)$,
``type'' is given by $(1/|1-c|, 1/|c-a-b|, 1/|a-b|)$.

\begin{table}[h!t!]
\caption{Special values of $_2F_1$ derived from closed-form relations for $F_1$.}
\label{value1}
\scalebox{1.00}{
\begin{tabular}
{llll}
\hline
No. & Formula & Type & Ref. \\ 
\hline
(A$''$.1)
& 
\shortstack{{} \\
${\displaystyle
\pFq{2}{1}{1/2,3/4}{1}{\frac {1}{81}}=
\frac{9}{100}
\frac{\sqrt{2}\, 5^{3/4}\,\G(1/4)^2}{\pi ^{3/2}}.
}$}
& (4,4,$\infty$) \\
(A$''$.2)
& 
\shortstack{{} \\
${\displaystyle
\pFq{2}{1}{1/2,3/4}{5/4}{\frac {80}{81}}=
\frac{9}{200}
\frac{5^{3/4}\,\G(1/4)^4}{\pi ^{2}}.
}$}
& (4,4,$\infty$) 
& *1
\\
(A$''$.3)
& 
\shortstack{{} \\
${\displaystyle
\pFq{2}{1}{1/4,1/2}{3/4}{\frac{80}{81}}=\frac{9}{5}.
}$}
& (4,4,$\infty$) 
& *2
\\
(B$''$.1)
& 
\shortstack{{} \\
${\displaystyle
\pFq{2}{1}{1/3,2/3}{7/6}{\frac {5}{32}}=
\frac{1}{20}\frac{\sqrt{3}\, \G(1/3)^6}{\pi ^3}.
}$}
& (3,6,6) 
\\
(B$''$.2)
& 
\shortstack{{} \\
${\displaystyle
\pFq{2}{1}{1/6,1/2}{5/6}{\frac {5}{32}}=\frac{2}{5}
\sqrt [6]{2}\sqrt {3}\sqrt [6]{5}. 
}$}
& (3,6,6) 
& *3
\\
(B$''$.3)
& 
\shortstack{{} \\
${\displaystyle
\pFq{2}{1}{1/3,2/3}{5/6}{\frac {27}{32}}
=\frac{8}{5}.
}$}
& (3,6,6) 
\\
(B$''$.4)
& 
\shortstack{{} \\
${\displaystyle
\pFq{2}{1}{1/2,5/6}{7/6}{\frac {27}{32}}=
\frac {1}{20}
\frac{\sqrt[3]{2}\sqrt{6}\,\G(1/3)^6}{\pi ^3}.
}$}
& (3,6,6) 
\\
(C$''$.1)
& 
\shortstack{{} \\
${\displaystyle
\pFq{2}{1}{3/5,4/5}{9/10}{\frac {3}{128}}
={\frac {8}{75}}\sqrt {2}\sqrt[10]{3}\sqrt {5}\sqrt {5+\sqrt {5}}.
}$}
& (2,5,10) 
\\
(C$''$.2)
& 
\shortstack{{} \\
${\displaystyle
\pFq{2}{1}{7/10, 9/10}{11/10}{\frac {3}{128}}
={\frac {8}{75}}\, 
\frac {\sqrt {2}\, \G(1/5)^3\,\G(2/5)}{(\sqrt{5}+1)\pi ^2}.
}$}
& (2,5,10) 
\\
(C$''$.3)
& 
\shortstack{{} \\
${\displaystyle
\pFq{2}{1}{3/5,4/5}{3/2}{\frac {125}{128}}=
\frac{4}{3}\,
\frac{\sqrt{2}\sqrt[10]{3}\,\G(1/5)^3\, \G(2/5)}{(5+\sqrt{5})^{3/2}\pi ^2}.
}$}
& (2,5,10) 
\\
(C$''$.4)
& 
\shortstack{{} \\
${\displaystyle
\pFq{2}{1}{1/10,3/10}{1/2}{\frac {125}{128}}=\frac{1}{6}\,
{2}^{3/10}\sqrt[10]{3}\left( \sqrt {5}+3 \right).
}$}
& (2,5,10) 
\\
(D$''$.1)
& 
\shortstack{{} \\
${\displaystyle
\pFq{2}{1}{1/6,2/3}{1}{\frac {9}{25}}
=\frac{1}{4}\,
\frac{\sqrt[3]{2}\sqrt[3]{5}\G(1/3)^3}{\pi ^2}.
}$}
& (2,6,$\infty$)
& *4 
\\
(D$''$.2)
& 
\shortstack{{} \\
${\displaystyle
\pFq{2}{1}{1/6, 2/3}{5/6}{\frac {16}{25}}=\frac{2}{3}\,
\sqrt [3]{5}.
}$}
& (2,6,$\infty$) 
& *5
\\
(D$''$.3)
& 
\shortstack{{} \\
${\displaystyle
\pFq{2}{1}{1/3,5/6}{7/6}{\frac {16}{25}}=
\frac{1}{48}\, \frac{\sqrt{3}\,5^{2/3}\,\G(1/3)^6}{\pi ^3}.
}$}
& (2,6,$\infty$)
& *6
\\
(E$''$.1)
& 
\shortstack{{} \\
${\displaystyle
\pFq{2}{1}{1/3,1/2}{1}{\frac {1}{5}}=
\frac{3}{20}\,
\frac{2^{2/3}\sqrt{5}\,\G(1/3)^3}{\pi ^2}.
}$}
& (6,6,$\infty$)\\
\hline
\end{tabular}
}
\scalebox{1}{
\begin{tabularx}{\linewidth}{X}
\small
*1:\, (7.10) in \cite{JZ2}.
\\
\small
*2:\, (1.6) in \cite{JZ2}.
\\
\small
*3:\, See *3 in Table {\ref{value2}} and *1 in Table {\ref{value3}}.
\\
\small
*4:\, See *4 in Table {\ref{value2}}.
\\
\small
*5:\, See *4 and *5 in Table {\ref{value2}}.
\\
\small
*6:\, See *4 and *6 in Table {\ref{value2}}.
\\
\hline
\end{tabularx}
}
\end{table}

\clearpage

\section{Algebraic transformations of identities in Table \ref{value1}}
In this section, applying algebraic transformations to the identities in Table \ref{value1},
we have some complicated identities for $_2F_1$.
These identities are tabulated in Tables \ref{value2} and \ref{value3}.

As an example, here we derive (A$'''$.1) in Table \ref{value2}.
As seen in Table \ref{value1}, (A$''$.1) is known.
Then, applying Goursat's algebraic transformation
\begin{gather*}
\pFq{2}{1}{a,b}{2b}{x}
=
(1-x)^{b-a}\left(1-\frac{x}{2}\right)^{a-2b}
\pFq{2}{1}{b-a/2,b+1/2-a/2}{b+1/2}{\left(\frac{x}{2-x}\right)^2}
\end{gather*}
(see formula (45) in \cite{Gour}) with $(a,b)=(3/4,1/2)$
to the left-hand side of (A$''$.1), 
(A$'''$.1) is easily deduced.
In a similar way, 
we derive the following:
\begin{itemize}
\item (A$'''$.2) from (A$''$.2) using
 Goursat's algebraic transformation (51) in \cite{Gour},
\item (A$'''$.3) from (A$''$.3) using
 Goursat's algebraic transformation (51) in \cite{Gour},
\item (B$'''$.1) from (B$''$.1) using
 Goursat's algebraic transformation (41) in \cite{Gour},
\item (B$'''$.2) from (B$''$.2) using
 Goursat's algebraic transformation (38) in \cite{Gour},
\item (B$'''$.3) using a connection formula for $_2F_1$ and (B$'''$.1) and (B$'''$.2),
\item (B$'''$.4) using a connection formula for $_2F_1$ and (B$'''$.1) and (B$'''$.2),
\item (C$'''$.1) from (C$''$.1) using
 Goursat's algebraic transformation (121) in \cite{Gour},
\item (C$'''$.2) from (C$''$.2) using
 Goursat's algebraic transformation (121) in \cite{Gour},
\item (C$'''$.3) using a connection formula for $_2F_1$ and (C$'''$.1) and (C$'''$.2),
\item (C$'''$.4) using a connection formula for $_2F_1$ and (C$'''$.1) and (C$'''$.2),
\item (D$'''$.1) from (D$''$.1) using
 Goursat's algebraic transformation (50) in \cite{Gour},
\item (D$'''$.2) from (D$''$.2) using
 Goursat's algebraic transformation (45) in \cite{Gour},
\item (D$'''$.3) from (D$''$.3) using
 Goursat's algebraic transformation (45) in \cite{Gour},
\item (E$'''$.1) from (E$''$.1) using
 Goursat's algebraic transformation (45) in \cite{Gour}.
\end{itemize}
The identities in Table \ref{value3} are derived as follows:
\begin{itemize}
\item (B$''''$.1) from (B$'''$.1) using
 Goursat's algebraic transformation (50) in \cite{Gour},
\item (B$''''$.2) from (B$'''$.2) using
 Goursat's algebraic transformation (50) in \cite{Gour},
\item (B$''''$.3) from (B$'''$.3) using
 Goursat's algebraic transformation (45) in \cite{Gour},
\item (B$''''$.4) from (B$'''$.4) using
 Goursat's algebraic transformation (45) in \cite{Gour}.
\end{itemize}

\begin{table}[h!t!]
\caption{Algebraic transformations of identities in Table \ref{value1}.} 
\label{value2}
\scalebox{1.00}{
\begin{tabular}{llll}
\hline
No. & Formula & Type & Ref. \\ 
\hline
(A$'''$.1)
& 
\shortstack{{} \\
$
{\displaystyle
\pFq{2}{1}{1/8, 5/8}{1}{\frac {1}{25921}}=\frac{1}{10}\,
\frac {\sqrt [4]{2}\sqrt [4]{161}\, \G(1/4)^2}{\pi ^{3/2}}.
}
$}
& (2,4,$\infty$) \\
(A$'''$.2)
&
\shortstack{{} \\ 
$
{\displaystyle
\pFq{2}{1}{3/8,7/8}{5/4}{\frac {25920}{25921}}
={\frac {1}{600}}\,\frac {{5}^{3/4}{161}^{3/4}\,\G(1/4)^4}{\pi ^2}.
}
$}
& (2,4,$\infty$) 
& *1
\\
(A$'''$.3)
&
\shortstack{{} \\ 
$
{\displaystyle
\pFq{2}{1}{1/8,5/8}{3/4}{\frac {25920}{25921}}=
\frac{3}{5}\,\sqrt [4]{161}. 
}
$}
& (2,4,$\infty$) 
&
*2
\\
(B$'''$.1)
& 
\shortstack{{} \\
$
{\displaystyle
\pFq{2}{1}{1/4,5/12}{7/6}{\frac {135}{256}}
={\frac {1}{20}}\,\frac {{2}^{5/6}\G(1/3)^6}{\pi ^3}.
}
$}
& (2,6,6) 
\\
(B$'''$.2)
& 
\shortstack{{} \\
$
{\displaystyle
\pFq{2}{1}{1/12,1/4}{5/6}{\frac {135}{256}}=\frac{2}{5}\,
\,\sqrt [6]{2}\sqrt {3}\sqrt [6]{5}.
}
$}
& (2,6,6) 
& *3
\\
(B$'''$.3)
&
\shortstack{{} \\ 
$
{\displaystyle
\pFq{2}{1}{1/4,5/12}{1/2}{\frac{121}{256}}=
\frac{1}{5}\,\frac{2^{5/12}\sqrt[8]{3}(1+\sqrt{3})^{3/2}\,\G(1/3)^3}{\sqrt{\pi}\, 
\G(1/4)^2}.
}
$}
& (2,6,6) 
\\
(B$'''$.4)
& 
\shortstack{{} \\
$
{\displaystyle
\pFq{2}{1}{3/4,11/12}{3/2}{\frac{121}{256}}=
\frac{8}{55}\,\frac{2^{11/12}\,3^{3/8}\, \G(1/3)^3\G(1/4)^2}
{(1+\sqrt{3})^{3/2}\pi^{5/2}}.
}
$}
& (2,6,6) 
\\
(C$'''$.1)
& 
\shortstack{{} \\
$
{\displaystyle
\pFq{2}{1}{1/30, 8/15}{9/10}{\frac {20736}{45125}}
=\frac{1}{6}\,
\sqrt {2}\sqrt[10]{3}\sqrt[10]{5}\sqrt[15]{19}\sqrt {5+\sqrt {5}}.
}
$}
& (2,3,10) 
\\
(C$'''$.2)
& 
\shortstack{{} \\
$
{\displaystyle
\pFq{2}{1}{2/15, 19/30}{11/10}{\frac {20736}{45125}}
=\frac{1}{12}\,
\frac{5^{2/5}19^{4/15}\G(1/5)^3\, \G(2/5)}
{(5+\sqrt{5})\, \pi ^2}.
}
$}
& (2,3,10) 
\\
(C$'''$.3)
& 
\shortstack{{} \\
$
{\displaystyle
\pFq{2}{1}{1/30,8/15}{2/3}{\frac{24389}{45125}}
=\frac{1}{36}\,
\frac{\sqrt{2}\sqrt{3}\sqrt{5+\sqrt{5}}+3\sqrt{5}+1}{\sqrt{5+\sqrt{5}}}
}
$}
& (2,3,10) 
\\
& 
\shortstack{{} \\
$
{\displaystyle
\phantom{\pFq{2}{1}{1/30,8/15}{2/3}{\frac{24389}{45125}}}
\times
\frac{2^{17/30}3^{9/10}5^{3/5}\sqrt[15]{19}\, \G(1/5)^2}
{\G(1/3)\G(1/15)}.
}
$}
\\
(C$'''$.4)
& 
\shortstack{{} \\
$
{\displaystyle
\pFq{2}{1}{11/30,13/15}{4/3}{\frac{24389}{45125}}
=\frac{5}{522}\,
\frac{\sqrt{5}-\sqrt{15-6\sqrt{5}}}{\sqrt{5}+\sqrt{15-6\sqrt{5}}}
}
$}
& (2,3,10) 
\\
&
\shortstack{{} \\
$
{\displaystyle
\phantom{\pFq{2}{1}{11/30,13/15}{4/3}{\frac{24389}{45125}}}
\times \frac{2^{7/30}\sqrt[10]{3}\,5^{3/5}19^{11/15}}
{\sqrt{5+\sqrt{5}}}
}
$}
\\
&
\shortstack{{} \\
$
{\displaystyle
\phantom{\pFq{2}{1}{11/30,13/15}{4/3}{\frac{24389}{45125}}}
\times
\frac{\G(1/5)\G(2/5)\G(1/3)\G(1/15)}{\pi^2}.
}
$}
\\
(D$'''$.1)
& 
\shortstack{{} \\
$
{\displaystyle
\pFq{2}{1}{1/3,1/3}{1}{\frac{1}{9}}=\frac{1}{4}\,
\frac{3^{2/3}\,\G(1/3)^3}{\pi^2}.
}
$}
&
$(3, \infty, \infty)$
& *4
\\
(D$'''$.2)
& 
\shortstack{{} \\
$
{\displaystyle
\pFq{2}{1}{1/3,1/3}{2/3}{\frac{8}{9}}=\frac{2}{3}\,3^{2/3}.
}
$}
&
$(3, \infty, \infty)$
& *5
\\
(D$'''$.3)
& 
\shortstack{{} \\
$
{\displaystyle
\pFq{2}{1}{2/3,2/3}{4/3}{\frac{8}{9}}=\frac{1}{16}\,
\frac{3^{5/6}\,\G(1/3)^6}{\pi^3}.
}
$}
&
$(3, \infty, \infty)$
& *6
\\
(E$'''$.1)
& 
\shortstack{{} \\
$
{\displaystyle
\pFq{2}{1}{1/6,2/3}{1}{\frac{1}{81}}=\frac{3}{20}\,
\frac{\sqrt[3]{2}\, 3^{2/3} \sqrt[6]{5}\, \G(1/3)^3}{\pi ^2}.
}
$}
&
$(2, 6, \infty)$\\
\hline
\end{tabular}
}
\scalebox{1}{
\begin{tabularx}{\linewidth}{X}
\small
*1:\, Formula 171 in Section 8.1.1 of \cite{Br}.
\\
\small
*2:\, (7.10) in \cite{JZ3}.
\\
\small
*3:\, See *1 in Table {\ref{value3}}.
\\
\small
*4:\, This is a special case of (1,2,3-1)(ix) in \cite{Eb}, 
(1/9.4) in \cite{Gos} and (1.2) in \cite{Ka}.
Because (D$''$.1) can be obtained
by applying an algebraic transformation 
to the left-hand side of (D$'''$.1),
it is regarded as a previously known formula.
The same holds for (D$''$.2) and (D$''$.3).
\\
\small
*5:\, This is a special case of (1,2,3-1)(vii) in \cite{Eb},
(8/9.1) in \cite{Gos} and (3.2) in \cite{Ka}.
\\
\small
*6:\, This is a special case of (1,2,3-1)(xxiii) in \cite{Eb},
(8/9.2) in \cite{Gos}, (5.24) in [DS], (3.3) in \cite{Ka}.
\\
\hline
\end{tabularx}
}
\end{table}

\clearpage

\begin{table}[!ht]
\caption{Algebraic transformations of identities in Table \ref{value2}.} 
\label{value3}
\scalebox{1.00}{
\begin{tabular}{llll}
\hline
No. & Formula & Type & Ref. \\ 
\hline
(B$''''$.1)
&
\shortstack{{} \\
$
{\displaystyle
\pFq{2}{1}{5/24,17/24}{7/6}{\frac {138240}{152881}}=\frac{1}{160}\,
\frac{\sqrt{2}\,17^{5/12}23^{5/12}\G(1/3)^6}{\pi ^3}.
}
$}
&(2.4.6)
\\
(B$''''$.2)
&
\shortstack{{} \\
$
{\displaystyle
\pFq{2}{1}{1/24, 13/24}{5/6}{\frac{138240}{152881}}=\frac{1}{5}\,
\sqrt{2}\sqrt{3}\sqrt[6]{5}\sqrt[12]{17}\sqrt[12]{23}.
}
$}
&(2.4.6)
& *1
\\
(B$''''$.3)
&
\shortstack{{} \\
$
{\displaystyle
\pFq{2}{1}{5/24,17/24}{3/4}{\frac {14641}{152881}}=\frac{1}{80}\,
2^{2/3}\sqrt[8]{3}\, 17^{5/12}23^{5/12}(1+\sqrt{3})^{3/2}
}
$}
&(2.4.6)
\\
&
\shortstack{{} \\
$
\displaystyle{
\phantom{
\pFq{2}{1}{5/24,17/24}{3/4}{\frac {14641}{152881}}}
\times
\frac{\G(1/3)^3}
{\sqrt{\pi}\,\G(1/4)^2}.}
$}
\\
(B$''''$.4)
&
\shortstack{{} \\
$
{\displaystyle
\pFq{2}{1}{11/24,23/24}{5/4}{\frac{14641}{152881}}=\frac{1}{1760}\,
\frac{2^{2/3}3^{3/8}17^{11/12}23^{11/12}}{(1+\sqrt{3})^{3/2}}
}
$}
&(2.4.6)
\\
&
\shortstack{{} \\
$
\displaystyle
\phantom{
\pFq{2}{1}{11/24,23/24}{5/4}{\frac{14641}{152881}}}
\times
\frac{\G(1/3)^3\G(1/4)^2}{\pi ^{5/2}}.
$}
\\
\hline
\end{tabular}
}
\scalebox{1}{
\begin{tabularx}{\linewidth}{X}
{\small
*1:\, We see that (B$''''$.2) is identical to (13) in \cite{BG} 
using the Pfaff transformation.
Therefore,
because (B$''$.2) and (B$'''$.2) can be obtained
by applying algebraic transformations to 
the left-hand side of (B$''''$.2),
these are regarded as previously known formulae.
}
\\
\hline
\end{tabularx}
}
\end{table}

\section{Concluding remarks}
Applying the method of contiguity relations to Appell's hypergeometric series $F_1$,
we have obtained several identities for $_2F_1$ with no free parameters.
In addition to the identities derived above,
there are a number of cases in which
the same method allows us to conjecture values of $_2F_1$ with no free parameters,
as we now discuss.

For example, for ${\bm k}=(-2, 3, 0, 1)$, 
we have the closed-form relation
\begin{gather*}
\frac{F(a+1)}{F(a)}=
-\frac{5}{3^3}\frac{(a+1/2)(a+3/2)}{(a+5/6)(a+7/6)},
\label{conj1}
\end{gather*}
where
\begin{gather*}
F(a):=\pFq{}{1}
{-2a; 3a+1. 1/2}{a+3/2}{\frac{16}{25}, \frac{4}{5}}.
\end{gather*}
From this relation, we can conjecture
\begin{align}
\begin{split}
F(a)&=\pFq{}{1}
{-2a; 3a+1. 1/2}{a+3/2}{\frac{16}{25}, \frac{4}{5}}\\
&=\left(\frac{5}{3^3}\right)^{a}
\frac{\cos (\pi a)\, \G(5/6)\G(7/6)\G(a+1/2)\G(a+3/2)}
{\G(1/2)\G(3/2)\G(a+5/6)\G(a+7/6)}\\
&=\frac{2}{3}\left(\frac{5}{3^3}\right)^{a}
\frac{\cos(\pi a)\,\G(a+1/2)\G(a+3/2)}{\G(a+5/6)\G(a+7/6)}.
\end{split}
\label{conj2}
\end{align}
Unfortunately, however, it seems difficult to
compute the asymptotic behavior of $F(a)$ in the limit $|a|\rightarrow +\infty$
from a direct application of Laplace's method.
For this reason, we have not been able to prove (\ref{conj2}).
However, 
with numerical calculations, we have obtained results consistent
with this identity.
If indeed this identity does hold, then we have
\begin{gather}
\pFq{2}{1}{1/2, 2/3}{7/6}{\frac{4}{5}}
=
\frac{1}{40}\frac{\sqrt{3}\, 5^{2/3}\, \G(1/3)^6}{\pi^3}
\label{conj3}
\end{gather}
by substituting $a=-1/3$ into (\ref{conj2}).
Moreover, using the connection formula (\ref{connection}), (E$''$.1) and (\ref{conj3}),
we obtain
\begin{gather}
\pFq{2}{1}{1/3, 1/2}{5/6}{\frac{4}{5}}=\frac{3}{\sqrt{5}}.
\label{conj4}
\end{gather} 
Although this relation follows from a conjecture, if indeed it does hold,
it provides a new example of an algebraic value of $_2F_1$. 
In addition, by applying algebraic transformations to (\ref{conj3}) and (\ref{conj4}),
we can similarly conjecture the relations
\begin{gather}
\pFq{2}{1}{1/3,5/6}{7/6}{\frac{80}{81}}=\frac{3}{40}\,
\frac{3^{5/6}\, \G(1/3)^6}{\pi ^3}
\label{conj5}
\end{gather}
and
\begin{gather}
\pFq{2}{1}{1/6,2/3}{5/6}{\frac{80}{81}}=\frac{3}{5}\,
3^{2/3}\, \sqrt[6]{5},
\label{conj6}
\end{gather}
respectively.

We can also obtain non-trivial algebraic values of $F_1$ 
from our approach.
For instance, substituting $a=-2/3$ into (B$'$.4),
we have
\begin{gather}
\pFq{}{1}{1/2; 2/3, 1/2}{11/6}{\frac{27}{32}, \frac{5}{6}}
=
\frac{2}{3}\,\sqrt{2}\sqrt{3}.
\end{gather}

We thus conclude that for the purpose 
of finding algebraic values of Gauss's hypergeometric series $_2F_1$ 
and Appell's hypergeometric series $F_1$,
it is worthwhile studying $F_1$ possessing closed forms.

\section*{Acknowledgement}
The author would like to thank Professor Katsunori Iwasaki
for many valuable comments.
The author is also grateful to Professor Frits Beukers 
for his encouragement.
This work is supported by a Grant-in-Aid for JSPS Fellows, JSPS No. 15J00201.

\medskip
\begin{flushleft}
Akihito Ebisu\\
Department of Mathematics\\
Hokkaido University\\
Kita 10, Nishi 8, Kita-ku, Sapporo, 060-0810\\
Japan\\
a-ebisu@math.sci.hokudai.ac.jp
\end{flushleft}


\begin{thebibliography}{KR}
\bibitem[Ar]{Ar}
N.Archinard,
{\it Exceptional sets of hypergeometric series},
J. Number Theory {\bf 101}(2003), no. {\bf 2}, 244--269. 
\bibitem[Ba]{Ba} W.N.Bailey, {\it Generalized hypergeometric series},
Cambridge Mathematical Tract No. 32, Cambridge University Press, (1935).
\bibitem[BG]{BG}
S.Baba, H.Granath,
{\it Quaternionic modular forms and exceptional sets 
of hypergeometric functions}, 
Int. J. Number Theory {\bf 11}(2015), no. {\bf 2}, 631--643.
\bibitem[Br]{Br}
Y.A.Brychkov,
{\it Handbook of special functions. 
Derivatives, integrals, series and other formulas},
CRC Press, Boca Raton, FL, 2008. 
\bibitem[BW]{BW} F.Beukers  and J.Wolfart,
{\it Algebraic values of hypergeometric functions}
in {\it New advances in transcendence theory (Durham, 1986)}, 
Cambridge Univ. Press, Cambridge(1988), 68--81 . 
\bibitem[Eb]{Eb} A.Ebisu,
{\it{Special values of the hypergoemetric series}},
Mem. Amer. Math. Soc., (to appear).
also available e-Print arXiv:1308.5588.
\bibitem[Erd1]{Erd1}
A.Erd\'elyi (editor), {\it Higher transcendental functions Vol.1}, 
McGraw-Hill, (1953).
\bibitem[Erd2]{Erd2}
A.Erd\'elyi, {\it Asymptotic expansions}
Dover Publications, Inc., New York, 1956.
\bibitem[Gos]{Gos} R.W.Gosper, {\it A letter to D. Stanton}, XEROX Palo Alto Research Center, 21
December 1977.
\bibitem[Gour]{Gour} \'{E}.Goursat, 
{\it Sur l'\'{e}quation diff\'{e}rentielle lin\'{e}aire, 
qui admet pour int\'{e}grale la s\'{e}rie hyperg\'{e}om\'{e}trique}, 
Ann. Sci. \'{E}cole Norm. Sup. ({\bf{2}}) {\bf{10}}(1881), 3--142. 
\bibitem[GS]{Gs} I.Gessel and D.Stanton, {\it Strange evaluations of hypergeometric series},
SIAM J. Math. Anal., {\bf{13}}(1982), no. {\bf{2}}, 295--308.
\bibitem[Iw]{Iw} K.Iwasaki, {\it On Some Hypergeometric Summations},
e-Print arXiv:1408.5658.
\bibitem[JZ1]{JZ1}
G.S.Joyce, I.J.Zucker, 
{\it Special values of the hypergeometric series},
Math. Proc. Cambridge Philos. Soc. {\bf 109}(1991), no. {\bf 2}, 257--261. 
\bibitem[JZ2]{JZ2}
G.S.Joyce, I.J.Zucker, 
{\it Special values of the hypergeometric series II},
Math. Proc. Cambridge Philos. Soc. {\bf 131}(2001), no. {\bf 2}, 309--319. 
\bibitem[JZ3]{JZ3}
G.S.Joyce, I.J.Zucker, 
{\it Special values of the hypergeometric series III},
Math. Proc. Cambridge Philos. Soc. {\bf 133}(2002), no. {\bf 2}, 213--222. 
\bibitem[Ka]{Ka} P.W.Karlsson, 
{\it On two hypergeometric summation formulas conjectured by Gosper}, 
Simon Stevin  {\bf {60}}  (1986),  no. 4, 329--337. 
\bibitem[Ki]{Ki} T.Kimura,
{\it Hypergeometric functions of two variables},
Seminar Note in Math. Univ. of Tokyo, 1973.
\bibitem[Ko]{Ko} W.Koepf, {\it Hypergeometric summation 
---An algorithmic approach to summation and special function identities, Second edition
}, Universitext, Springer(2014).
\bibitem[PWZ]{PWZ} M.Petkov\v{s}ek, H.Wilf and D.Zeilberger, {\it A=B}, A.K.Peters, Wellesley(1996).
\bibitem[Si]{Si} 
C.L. Siegel, 
{\it {\"{U}}eber einige Anwendungen diophantischer Approximationen}, 
in: Ges. Werke, Bd 1,
Springer, Berlin, Heidelberg, New York, 1966, 209--266.
\bibitem[Ta1]{Ta1} N.Takayama,
{\it Gr\"obner basis and the problem of contiguous relations}, 
Japan J. Appl. Math. {\bf 6}(1989), no. {\bf 1}, 147--160. 
\bibitem[Ta2]{Ta2} N.Takayama, 
{\it  Gr\"obner basis for rings of differential operators and applications},
in {\it  Gr\"obner bases: Statistics and software systems}, Springer(2013), 279--344. 
\bibitem[Vi]{Vi}
R.Vidunas,
{\it Expressions for values of the gamma function},
Kyushu J. Math. {\bf 59}(2005), no. {\bf 2}, 267--283.
\bibitem[Wo]{Wo}
J.Wolfart, 
{\it Werte hypergeometrischer Funktionen}, 
Invent. Math. {\bf 92}(1988), no. {\bf 1}, 187--216.  
\end{thebibliography}
\end{document}